\documentclass[11pt]{article}

\usepackage{amsmath}
\usepackage{amssymb}  
\usepackage{epsfig}
\usepackage{amsthm}   
\usepackage[mathcal]{eucal}
\newcommand{\F}{\mathcal{F}}

\newcommand{\R}{\mathbb{R}}
\newcommand{\BR}{\bar{\mathbb{R}}}

\DeclareMathOperator{\clconv}{cl\,conv}

\newcommand{\rot}{\Lambda}

\newcommand{\inner}[2]{\langle{#1},{#2}\rangle}

\newcommand{\norm}[1]{\|#1\|}

\newcommand{\tos}{\rightrightarrows} 

\newtheorem{theorem}{Theorem}[section]
\newtheorem{lemma}[theorem]{Lemma}

\newtheorem{definition}{Definition}[section]

\title{Maximal monotone operators with a unique extension to the bidual}

\author{M. Marques Alves\thanks{IMPA, Estrada Dona Castorina 110, 22460-320
    Rio de Janeiro, Brazil
   ({\tt maicon@impa.br})}\hspace{.5em}\thanks{Partially supported by Brazilian CNPq
    scholarship 140525/2005-0.}
  \and
    B. F. Svaiter\thanks{ IMPA, Estrada Dona Castorina 110, 22460-320 Rio de
    Janeiro, Brazil ({\tt benar@impa.br}) }\hspace{.5em}
    \thanks{Partially supported by CNPq
    grants 300755/2005-8, 475647/2006-8 and by PRONEX-Optimization}
}

\date{}

\begin{document}

\maketitle

\begin{abstract}
  We present a new sufficient condition under which a maximal monotone
  operator $T:X\tos X^*$ admits a unique maximal monotone
  extension to the bidual $\widetilde T:X^{**}\tos X^*$.
  For non-linear operators this condition is equivalent to uniqueness of
  the extension.

  The class of maximal monotone operators which satisfy this new condition
  includes class of Gossez type D maximal monotone operators, previously
  defined and studied by J.-P. Gossez,
  and all maximal monotone operators of this new class satisfies a
  restricted version of Br\o ndsted-Rockafellar condition.

  The central tool in our approach is the $\mathcal{S}$-function
  defined and studied by Burachik and Svaiter in
  2000~\cite{BuSvSet02}(submission date, July 2000).
  For a generic operator, this function is the supremum of all convex
  lower semicontinuous functions which are majorized by the duality
  product in the graph of the operator.

  We also prove in this work that if the graph of a maximal monotone
  operator is convex, then this graph is an affine linear subspace.
  \\
  \\
  2000 Mathematics Subject Classification: 47H05, 49J52, 47N10.
  \\
  \\
  Key words: Maximal monotone operators, extension,  bidual, Banach
  spaces, 
  Br\o ndsted-Rockafellar property, $\mathcal{S}$-function.
  \\
\end{abstract}

\pagestyle{plain}


\section{Introduction}

Let $X$ be a real Banach space.  We use the notation $X^*$ for the
topological dual of $X$ and $\pi_{X\times X^*}$,
$\inner{\cdot}{\cdot}_{X\times X^*}$ for the duality product
\[
\pi_{X\times X^*}(x,x^*)=\inner{x}{x^*}_{X\times X^*}=x^*(x).
\]
Whenever the underlying domain of the duality product is clear, we will use the
notations $\pi$ and $\inner{\cdot}{\cdot}$. 

A point  point-to-set operator $T:X\tos X^*$ 
(respectively $T:X^{**}\tos X^{*}$)
 is a relation on $X$ to
$X^*$ 
(respectively on $X^{**}$ to $X^{*}$):
\[ 
 T\subset X\times X^* \quad\mbox{(respectively } T\subset
X^{**}\times X^*),
\]
and $r\in T(q)$ means $(q,r)\in T$.
An operator $T:X\tos X^*$, or $T:X^{**}\tos X^*$, is {\it monotone} if
\[
\inner{q-q'}{r-r'}\geq 0,\forall (q,r),(q',r')\in T.
\]
An operator $T:X\tos X^*$ is {\it maximal monotone} (in $X\times X^*$) if it
is monotone and maximal (whit respect to the inclusion) in the family
of monotone operators of $X$ in $X^*$.
An operator $T:X^{**}\tos X^*$ is maximal monotone (in $X^{**}\times
X^*$) if it is monotone and maximal (with respect to the inclusion) in
the family of monotone operators of $X^{**}$ in $X^*$.

The canonical injection of $X$ into $X^{**}$ allows one to identity $X$
with a subset of $X^{**}$.
Therefore,
any maximal
monotone operator $T:X\tos X^*$ is also a monotone
operator $T:X^{**}\tos X^*$ and admits one (ore more)
maximal monotone extension in  $X^{**}\times X^*$.
In general this maximal monotone extension will
not be unique. We are concerned with the problem:
\begin{quotation}
  \noindent
  Under which conditions a maximal monotone operator  $T:X\tos X^*$ has
  a unique extension to the bidual,  $X^{**}\tos  X^*$?
\end{quotation}
The problem of unicity of maximal extension of a generic monotone
operator was studied in details by Legaz and Svaiter
in~\cite{legaz-svaiter-01}. That paper will be an important reference
for the present work.

The specific problem above mentioned, of uniqueness of extension a
\emph{maximal monotone} operator to the bidual, has been previously
addressed by
Gossez~\cite{GosProc71,GosProc72,GosConvProc76,GosExProc76}. He found
a condition under which uniqueness of the extension is
guaranteed~\cite{GosExProc76}.
We will provide another condition, which encompass Gossez's condition,
and we will also prove that operators which satisfy this more general
condition satisfy also a relaxed version of the Br\o ndsted-Rockafellar
condition~\cite{MASvJCA08}.

\begin{definition}
  \label{def:g}
A maximal monotone operator $T:X\tos X^*$ is \emph{Gossez type D} if for any
$(x^*,x^{**})\in X^*\times X^{**}$ such that
\[
 \inner{y-x^{**}}{y^*-x^*}\geq 0,\, \forall (y,y^*)\in T
\]
there exists a \emph{bounded} net $\{ (x_i,x_i^*)\}_{i\in I}\subset T$ such that
\[ x_i \stackrel{\sigma(X^{**},X^*)}{\longrightarrow} x^{**},\qquad
   x_i^*\stackrel{\norm{\cdot}}{\longrightarrow} x^*.
\]
\end{definition}
Gossez has proved that if $T\subset X\times X^{*}$ is maximal monotone of
Gossez type $ D$, then it has a \emph{unique} maximal monotone
extension to $X^{**}\times
X^*$, that is, $\widetilde T:X^{**}\tos X^*$ maximal monotone and $T\subset
\widetilde T$.

Now we will discuss the Br\o ndsted-Rockafellar property.  Let $T\subset
X\times X^*$ be maximal monotone.  Burachik, Iusem and Svaiter~\cite{BuIuSv}
 defined
the $T^\varepsilon$ enlargement of $T$ for $\varepsilon\geq 0$, as
$T^\varepsilon :X\tos X^*$
\begin{equation}
  \label{eq:def.te}
   T^\varepsilon(x)=\{x^*\in X^*\;|\;  \inner{x-y}{x^*-y^*}\geq -\varepsilon
  \quad \forall (y,y^*)\in T\}.
\end{equation}
It is trivial to verify that $T\subset T^\varepsilon$. The
$T^\varepsilon$ enlargement is a generalization of the
$\varepsilon$-subdifferential. As the $\varepsilon$-subdifferential,
the $T^\varepsilon$ has also practical 
uses~\cite{SolodovSv1,SolodovSv2,Correa1,Str1,K1,Moud01}.
Br\o ndsted and Rockafellar proved that the the
$\varepsilon$-subdifferential may be seen as an approximation of the
\emph{exact} subdifferential at a nearby point. This property, may be
extended to the context of maximal monotone operators.  A maximal
monotone operator has the \emph{Br\o ndsted-Rockafellar property} if,
for any $\varepsilon,\lambda>0$,
\[
 x^*\in T^\varepsilon(x)\Rightarrow 
\forall \lambda>0,\;
\exists (\bar x, \bar x ^*)\in
T,\quad \norm{x-\bar x}\leq \lambda,\quad \norm{\bar
  x^*-x^*}\leq \varepsilon/\lambda.
\]
This property was proved in Hilbert spaces in~\cite{BuSaSv1}. Latter
on, Burachik and Svaiter~\cite{BuSvteb} proved that in a
\emph{reflexive} Banach space, all maximal monotone operators
satisfies this property.
The operator $T$ satisfy the \emph{restricted Br\o ndsted-Rockafellar
  property}~\cite{MASvJCA08} if 
\begin{align}
  \nonumber
x^*\in T^\varepsilon(x),\;\tilde \varepsilon>\varepsilon\Rightarrow
\forall \lambda>0,
&\exists (\bar x, \bar x ^*)\in T,\\
& \norm{x-\bar x}< \lambda,\quad \norm{\bar
  x^*-x^*}< \tilde \varepsilon/\lambda.
  \label{eq:r.br}
\end{align}
In a recent work~\cite{MASvJCA08} the authors defined a general class of maximal monotone operators in non-reflexive Banach spaces which satisfies the above property.
An open question is whether a maximal monotone operator Gossez type D
has Br\o ndsted-Rockafellar property.

We will analyze a new condition for, which embraces Gossez type D
condition,  and which is sufficient for a maximal monotone operator to:
\begin{enumerate}
\item admit a
unique maximal monotone extension to the bidual space;
\item satisfy the restricted Br\o ndsted-Rockafellar property.
\end{enumerate}
Moreover, we will prove that, for non-linear maximal monotone
operators,  this new  condition is equivalent to unicity of extension to
the bidual.

We will also study the relation of this new condition with general
properties of Fitzpatrick family (defined bellow). So, the main
properties of this families, systematically studied in the series~\cite{SvSet00,
  BuSvSet02,SvFixProc03,BuSvProc03,MASvJCA08} will also be used. As
this field has been subject of intense research, to clarify precedence
issues, we will include submission date of some works where these
properties were obtained.

Given a maximal monotone operator $T:X\tos X^*$, Fitzpatrick
defined~\cite{Fitz88} the family $\mathcal{F}_T$ as those convex, lower
semicontinuous functions in $X\times X^*$ which are bounded bellow by
the duality product and coincides with it at $T$:
\begin{equation} \label{eq:def.ft}
  \F_T=\left\{ h\in \BR^{X\times X^*}
    \left|
      \begin{array}{ll}
        h\mbox{ is convex and lower semicontinuous}\\
        \inner{x}{x^*}\leq h(x,x^*),\quad \forall (x,x^*)\in X\times X^*\\
        (x,x^*)\in T 
        \Rightarrow 
        h(x,x^*) = \inner{x}{x^*}
      \end{array}
    \right.
  \right\}.
\end{equation}
Fitzpatrick found an explicit formula for the minimal element of
$\mathcal{F}_T$, from  now on \emph{Fitzpatrick function} of $T$:
\begin{equation}
  \label{eq:def.f.fitz}
  \varphi_T(x,x^*)=\sup_{(y,y^*)\in T} \inner{x}{y^*}+\inner{y}{x^*}-
    \inner{y^*}{y}.
\end{equation}
Note that in the above definition, $T$ may be a generic subset of
$X\times X^*$.  Legaz and Svaiter in~\cite{legaz-svaiter-01} studied
generic properties of $\varphi_T$ for arbitrary sets and its
relation with monotonicity.

The \emph{conjugate} of a function $f:X\to\BR$ is
defined as $f^*:X^*\to\BR$,
\[
 f^*(x^*)=\sup_{x\in X} \inner{x}{x^*}-f(x),
\]
and the \emph{convex closure} of $f$ is $\clconv f:X\to\BR$,
the largest convex lower semicontinuous function majorized by $f$:
\[ 
 \clconv f(x):=\sup\{ h(x)\;|\;h \mbox{ convex, lower semicontinuous, }
  h\leq f\}.
\]
The indicator function of $A\subset X$ is $\delta_{A,X}:X\to\BR$,
\[\delta_{A,X}(x)=
 \begin{cases}
   0& x\in A\\
   \infty& \mbox{otherwise.}
 \end{cases}
\]
Whenever the set $X$ is implicitly defined, we use the notation
$\delta_A$.

\begin{definition}[\cite{BuSvSet02}(submission date, July
  2000, Eq.\ (35)]
  \label{def:s}\hspace{1em}\\
  The $\mathcal{S}$-function (original notation
  $\Lambda_{\overline{S^T}}$) associated with a maximal monotone
  operator $T:X\tos X^*$ is  $ \mathcal{S}_{\,T}:X\times X^*\to\BR$
  \begin{equation}
    \label{eq:def.f.s}
    {\mathcal{S}}_{\,T}=\clconv\left(\pi+\delta_T\right). 
  \end{equation}
\end{definition}
In \cite{BuSvSet02} it is proved that the $\mathcal{S}$-function is the
\emph{supremum} of the family of Fitzpatrick function.
The epigraphical structure of the $\mathcal{S}$-function was previously
studied in~\cite{SvSet00} (submission date, September 1999).
This function
will be central for defining the new class of maximal monotone
operators that admits a unique maximal monotone extension to the
bidual. 

The $\mathcal{S}$-function and Fitzpatrick function are still well defined for 
arbitrary sets (or operators) $T\subset X\times X^*$:
\begin{align}
  \label{eq:new.1}
     \mathcal{S}_T:X\times X^*\to\BR,\; &\;
     \mathcal{S}_T=\clconv (\pi+\delta_T),\\[1em]
  \label{eq:def.f.fitz.2}
     \varphi_T:X\times X^*\to\BR, &\;
       \varphi_T(x,x^*)=\sup_{(y,y^*)\in T} \inner{x}{y^*}+\inner{y}{x^*}-
    \inner{y^*}{y}.
\end{align}
Legaz and Svaiter also studied in~\cite{legaz-svaiter-01} generic
properties of $\mathcal{S}$ (with the notation $\sigma_T$) and $\varphi_T$ for
arbitrary sets and its relation with monotonicity and maximal
monotonicity.

 To simplify the notation, define
\[ \rot:X^{**}\times X^*\to X^*\times X^{**},\quad
\rot(x^{**},x^*)=(x^*,x^{**}).
\]
Note that $\rot(X\times X^*)=X^*\times X$.
We will prove three  main results in this paper:
\begin{theorem}
  \label{th:1}
  Let $X$ be a generic Banach space and $T:X\tos X^*$ a maximal
  monotone operator such that
 \begin{equation}
    \label{eq:cond.as}
    (\mathcal{S}_{T})^*(x^*,x^{**})\geq\inner{x^*}{x^{**}},\qquad \forall
    (x^*,x^{**})\in X^*\times X^{**}.
  \end{equation}
  Then $T$ admits a unique maximal monotone
  extension  $\widetilde T:X^{**}\tos X^*$.
  
  Additionally, 
  $(\mathcal{S}_{T})^*=\varphi_{\rot \widetilde T}$ and 
  for all $h\in \mathcal{F}_T$,
  \begin{align*}
  &   h^*(x^*,x^{**})\geq\inner{x^*}{x^{**}},\quad \forall
  (x^*,x^{**})\in X^*\times X^{**};\\
  &h^*\in \mathcal{F}_{\rot\,\widetilde{T}}\;.
  \end{align*}
  Moreover, $T$ satisfy the restricted Br\o ndsted-Rockafellar
  property.
\end{theorem}
The last statement of the above theorem is a particular case of a more
general result proved in~\cite{MASvJCA08}. In that paper, it is proved
that if a convex lower semicontinuous function in $X\times X^*$ and
its conjugate majorizes the duality product in $X\times X^*$ and
$X^{*}\times X^{**}$ respectively, then this function is in the
Fitzpatrick family of a maximal monotone operator and this maximal
monotone operator satisfy the restricted Br\o ndsted-Rockafellar
condition.  This result, in a reflexive Banach space was previously
obtained in~\cite{BuSvProc03}.

A natural question is whether the converse of Theorem~\ref{th:1}
holds. To give a partial answer to this question, first recall that
a linear (affine) subspace of a real linear space $Z$ is a set $A\subset Z$
such that there exists $V$, subspace of $Z$, and a point $z_0$ such that
\[
 A=V+\{z_0\}=\{z+z_0\;|\; z\in V\}
\]
\begin{theorem}
  \label{converse}
  Suppose that $T:X\tos X^*$ is maximal monotone and is not affine
  linear. In this case, if $T$ has a unique extension to $X^{**}\times
  X^{*}$, then
   \begin{equation}
    \label{eq:cond.as.same}
    (\mathcal{S}_T)^*(x^*,x^{**})\geq\inner{x^*}{x^{**}},\qquad \forall
    (x^*,x^{**})\in X^*\times X^{**},
  \end{equation}
  and $T$ satisfy the restricted Br\o ndsted-Rockafellar property.
\end{theorem}

According to the above theorems, for non-linear maximal monotone
operators, condition \eqref{eq:cond.as} is equivalent to unicity of
maximal monotone extension to the bidual.
Surprisingly, condition \eqref{eq:cond.as} is as general or even more
general that Gossez's type D property. This is the last result
of this work:
\begin{theorem}
  \label{th:main}
  Let $X$ be a generic real Banach space and $T:X\tos X^*$ be maximal
  monotone and Gossez type D. Then,
 \begin{equation*}
    (\mathcal{S}_T)^*(x^*,x^{**})\geq\inner{x^*}{x^{**}},\qquad \forall
    (x^*,x^{**})\in X^*\times X^{**}.
  \end{equation*}
 In particular, 
 for all $h\in \mathcal{F}_T$,
   \[ h^*(x^*,x^{**})\geq \inner{x^*}{x^{**}},\qquad \forall (x^*,x^{**})\in
     X^*\times X^{**}.
     \]
     and  $T$ satisfy the Br\o ndsted-Rockafellar property.
\end{theorem}
%
\section{Preliminary results} \label{sec:pre}

As mentioned before, Fitzpatrick proved that the family $\mathcal{F}_T$
is non-empty by producing its smallest element $\varphi_T$. 
Fitzpatrick also proved that \emph{any} function in this family fully
characterizes the maximal monotone operator which defines the family:
\begin{theorem}[\cite{Fitz88}]
  \label{th:fitz}
  Let $T:X\tos X^*$ be maximal monotone. Then, for any $h\in\mathcal{F}_T$,
  \[
  h(x,x^*)=\inner{x}{x^*}\iff (x,x^*)\in T,\qquad \forall (x,x^*)\in
  X\times X^*.
  \]
  Moreover, $\varphi_T$ is the smallest function of $\mathcal{F}_T$.
\end{theorem}

In~\cite{BuSvSet02}(submission date, July 2000), Burachik and Svaiter proved
\begin{theorem}[\mbox{\cite[eqns.\ (32), (37), (39) ]{BuSvSet02}}]
  \label{th:bs}
  Let $T:X\tos X^*$ be  maximal monotone  and $\mathcal{S}_T$ be the $\mathcal{S}$-function
  associated with $T$, as defined in \eqref{eq:def.f.s}. Then,
  \[ \mathcal{S}_T\in \mathcal{F}_T, \qquad  \mathcal{S}_{\,T}:=\sup_{h\in\F_T}\;\{h\}.
  \]
  and
  $\varphi_T\leq h\leq \mathcal{S}_T$ for all $h\in \mathcal{F}_T$.

   Moreover, $\varphi_T$ and $\mathcal{S}_T$ are related as follows:
  \[
  \varphi_T(x,x^*)=(\mathcal{S}_T)^*(x^*,x),\qquad \forall (x,x^*)\in X\times X^*.
  \]
\end{theorem}
Define, for $h:X\times X^*\to\BR$,
\begin{equation}
  \label{eq:def.jop}
  \mathcal{J}
     h:X\times X^*\to\BR,\quad  \mathcal{J}
     h (x,x^*)=h^*(x^*,x).
\end{equation}
According to the above theorem, $\mathcal{J} \mathcal{S}_T=\varphi_T\in
\mathcal{F}_T$.  So, it is natural to ask whether $\mathcal{J}$ maps
$\mathcal{F}_T$ in to itself. Burachik and Svaiter also proved that
this happens in fact:
\begin{theorem}[\mbox{\cite[Theorem 5.3]{BuSvSet02}}]

  \label{th:bs.2}
  Suppose that $T$ is maximal monotone. Then
  \[ \mathcal{J}h\in
 \mathcal{F}_T, \qquad \forall  h\in \mathcal{F}_T,\]
  that is, 
  if $h\in\mathcal{F}_T$, and
  \[ g:X\times X^*\to\BR,\quad g(x,x^*)=h^*(x^*,x),
  \]
  then $g\in\mathcal{F}_T$.

  In a reflexive Banach space
  $\mathcal{J}\varphi_T=\mathcal{S}_T$.
\end{theorem}
It is interesting to note that $\mathcal{J}$ is a order-reversing
mapping of $\mathcal{F}_T$ into itself. This fact suggests that this
mapping may have fixed points in $\mathcal{F}_T$. Svaiter proved
\cite{SvFixProc03}(submission date, July 2002) that if $T$ is maximal
monotone, then $\mathcal{J}$ has always a fixed point in
$\mathcal{F}_T$.

Legaz and Svaiter~\cite{legaz-svaiter-01} observed that for a generic $T\subset
X\times X^*$ 
\begin{equation}
  \label{eq:leg-sv01}
  \varphi_T(x,x^*)=(\pi+\delta_T)^*(x^*,x)=
  (\mathcal{S}_T)^*(x^*,x), \qquad \forall (x,x^*)\in X\times X^*,
\end{equation}
Therefore, also for an arbitrary $T$, $ \mathcal{J} \mathcal{S}_T=\varphi_T$.

It will be useful to consider monotonicity a \emph{relation} and to
study it also in the framework of the classical notion of
polarity~\cite{birk}
\begin{definition}[\cite{legaz-svaiter-01}]
  \label{def:polar}
  Two points $(x,x^*),(y,y^*)\in X\times X^*$ are in \emph{monotone
    relation}, $(x,x^*)\mu(y,y^*)$ if
  \[
  \inner{x-y}{x^*-y^*}\geq 0.
  \]
  Give $A\subset X\times X^*$, the \emph{monotone polar} (in $X\times
  X^*$) of $A$ is the set $A^{\mu}$,
  \begin{equation}
    \label{eq:def.mon.pl}
    \begin{array}{rcl}
      A^\mu&=&\{(x,x^*)\in X\times X^*\;|\; (x,x^*)\mu(y,y^*),\qquad 
      \forall (y,y^*)\in A\},\\
      &=&\{(x,x^*)\in X\times X^*\;|\; \inner{x-y}{x^*-y^*}\geq 0,\quad 
      \forall (y,y^*)\in A\}.
    \end{array}
  \end{equation} 
\end{definition}
We shall need some results of Legaz and Svaiter
which are scattered
along~\cite{legaz-svaiter-01} and which we expound in the next two theorems:
\begin{theorem}[\mbox{\cite[Eq. (22),
  Prop.\ 2,
  Prop.\ 21]{legaz-svaiter-01}}]
  \label{th:legaz-svaiter}
  Let $A\subset X\times X^*$. Then
  \begin{equation}
    \label{eq:22}
    A^\mu=\{(x,x^*)\in X\times X^*\;|\; \varphi_T(x,x^*)\leq
     \inner{x}{x^*}\},
   \end{equation}
   and  the following conditions are equivalent
 \begin{enumerate}
 \item $A$ is monotone,
 \item $\varphi_A\leq (\pi+\delta_A)$.
 \item $A\subset A^\mu$,
 \end{enumerate}
 \end{theorem}
 Note in the above theorem and in the definition of Fitzpatrick's
 family, the convenience of defining as in~\cite[Eq. (12) and
 bellow]{legaz-svaiter-01}, for $h:X\times X^*\to\BR$:
 \begin{equation}
   \label{eq:def.bL}
   \begin{array}{rl}
     b(h)&:=\{ (x,x^*)\in X\times X^*\;|\;
       h(x,x^*)\leq \inner{x}{x^*}\},\\
    L(h)&:=\{ (x,x^*)\in X\times X^*\;|\;
       h(x,x^*) = \inner{x}{x^*}\}.
   \end{array}
 \end{equation}
 \begin{theorem}[\mbox{\cite[Prop. 36,
  Lemma 38]{legaz-svaiter-01}}]
   \label{th:ls-02}
  Suppose that $A\subset X\times X^{*}$ is monotone.
  Then the following conditions are equivalent
  \begin{enumerate}
   \item $A$ has a unique maximal monotone extension (in $X\times X^*$),
   \item $A^\mu$ is monotone
   \item $A^\mu$ is maximal monotone,
  \end{enumerate}
  and  if any of these conditions holds, then $A^\mu$ is the
  unique maximal monotone extension of $A$.

  Moreover, still assuming only $A$ monotone,
  \begin{equation}
    \label{eq:a1}
    \varphi_A\geq \pi\iff
    b(\varphi_A)=L(\varphi_A)
  \end{equation}
 and if these conditions holds, then $A$ has a unique maximal
 monotone extension, $A^\mu$.
\end{theorem}


\section{Convexity and maximal monotonicity}
\label{sec:mmc}

To prove Theorem~\ref{converse} we shall need an additional result,
which, up to our knowledge, is also new:

\begin{lemma}
  \label{bas}
  If $T:X\tos X^*$ is maximal monotone and convex, then $T$ is affine linear.
\end{lemma}
\begin{proof}
   Take and arbitrary $(x_0,x^*_0)\in T$ and define 
   \[ T_0=T-\{(x_0,x^*_0)\}
   \]
   Note that $T_0$ is maximal monotone and convex.  So, it suffices to
   prove that $T_0$ is a linear subspace of $X\times X^*$. Take an
   arbitrary $(x,x^*)\in T_0$. First we claim that
   \begin{equation}
     \label{eq:claim.1}
     t(x,x^*)\in T_0,\qquad \forall  t\geq 0.
   \end{equation}
   For $0\leq t\leq 1$ the above inclusion holds because $(0,0)\in
   T_0$ and $T_0$ is convex
   For the case $t\geq 1$ let $(y,y^*)\in T$. Then,
   $t^{-1}(y, y^*)\in T_0$ and so
   \[
    \inner{x-t^{-1}y}{x^*-t^{-1}y^*}\geq 0.
    \]
    Multiplying this inequality by $t$ we conclude that
    \( \inner{tx-y}{tx^*-y^*}\geq 0\). As $(y,y^*)$ is a generic
    element of $T_0$, which is maximal monotone, we conclude that
    $t(x,x^*)\in T_0$ and the claim \eqref{eq:claim.1} holds.

    We have just proved that $T_0$ is a convex cone.
    Now take an
    arbitrary pair
    \[ (x,x^*),(y,y^*)\in T_0.\]
    Then
    \begin{equation}
      \label{eq:claim.2}
         (x+y,x^*+y^*)=
    2\left[\frac{1}{2}(x,x^*)+\frac{1}{2}(y,y^*)\right]\in T_0. 
    \end{equation}
    As $(0,0)\in T_0$,
    \[
    \inner{y-(-x)}{y^*-(-x^*)}=\inner{(y+x)-0}{(y^*+x^*)-0}\geq 0.
    \]
    Since $T_0$ is maximal monotone, we conclude that $-(x,x^*)\in
    T_0$.
    Therefore, using again~\eqref{eq:claim.1} we conclude that $T_0$
    is closed under scalar multiplication. To end the proof, combine
    this result with~\eqref{eq:claim.2} to conclude that $T_0$ is a
    linear subspace.
\end{proof}

\section{Proof of the main results}
\label{sec:proofs}

From now on, $T:X\tos X^*$ is a maximal monotone operator.  The inverse
of $T$ is $T^{-1}:X^*\tos X$,
\begin{equation}
  \label{eq:def.t-1}
  T^{-1}=\{ (x^*,x)\in X^*\times X\;|\; (x,x^*)\in T\}.
\end{equation}
Note that $T^{-1}\subset X^*\times X\subset X^{*}\times X^{**}$.
Fitzpatrick function of $T^{-1}$, regarded as a subset of $X^*\times X^{**}$
is, according to~\eqref{eq:def.f.fitz.2}
\begin{align*}
  \varphi_{T^{-1}, X^*\times X^{**}}(x^*,x^{**})&=
  \sup_{( y^*,y^{**})\in T^{-1}}\inner{x^*}{y^{**}}+\inner{y^*}{x^{**}}
  -\inner{y^*}{y^{**}}\\
  &=
  \sup_{( y^*, y)\in T^{-1}}\inner{x^*}{y}+\inner{y^*}{x^{**}}
  -\inner{y^*}{y}\\
  &=\left(\pi+\delta_T\right)^*(x^*,x^{**}).
\end{align*}
where the last $x^*$ is identifying with its image under the
canonical injection of $X^*$ into $X^{***}$. 
Using the above equations, \eqref{eq:new.1} and the fact that conjugation
is invariant under the convex-closure operation we obtain
  \begin{equation}
    \label{eq:ft.tm2}
     \varphi_{T^{-1}, X^*\times X^{**}}= \left(\pi+\delta_T\right)^*=
     (\mathcal{S}_T)^*.
   \end{equation}
   where $\pi=\pi_{X\times X^*}$ and $\delta_T=\delta_{T,X\times X^*}$.

   We will use the notation $\left(T^{-1}\right)^{\mu,X^*\times X^{**}}$
   for denoting the monotone polar of $T^{-1}$ in
   $X^*\times X^{**}$.
  Combing the above equation with Theorem~\ref{th:legaz-svaiter} we
  obtain a simple expression for this  monotone polar:
  \begin{align}
    \label{eq:w2}
    \left(T^{-1}\right)^{\mu,X^*\times X^{**}} &=\{(x^{*},x^{**})\in
    X^{*}\times X^{**}\,| \, (\mathcal{S}_T)^*(x^*,x^{**})\leq
    \inner{x^*}{x^{**}}\}.
  \end{align}

\begin{proof}[Proof of Theorem \ref{th:1}]
  Combining assumption~\eqref{eq:cond.as} and \eqref{eq:ft.tm2} we have
  \[
  \varphi_{T^{-1}, X^*\times
    X^{**}} (x^*,x^{**})=(\mathcal{S}_T)^*(x^*,x^{**})\geq\inner{x^*}{x^{**}},
  \qquad \forall (x^*,x^{**})\in X^{*}\times X^{**}.
  \]
  Therefore, using Theorem~\ref{th:ls-02} and
  Theorem~\ref{th:legaz-svaiter} for $A=T^{-1}\subset X^*\times X^{**}$ 
  we conclude that
  $\left(T^{-1}\right)^{\mu,X^*\times X^{**}}$, the monotone polar of
  $T^{-1}$ in $X^*\times X^{**}$, is the \emph{unique} maximal
  monotone extension of $T^{-1}$ to $X^*\times X^{**}$ and
  \begin{equation}
    \label{eq:t1.w}
    \left(T^{-1}\right)^{\mu,X^*\times X^{**}}
     =\{ (x^*,x^{**}) \in X^*\times X^{**}\,|\; (\mathcal{S}_T)^*(x^*,x^{**})
    =\inner{x^*}{x^{**}}\}.
  \end{equation}
  Using the above result and again \eqref{eq:cond.as}, we conclude
  that 
  \[  (\mathcal{S}_T)^*\in\mathcal{F}_{ \left(T^{-1}\right)^{\mu,X^*\times X^{**}}}.
  \]
  Now,   define
  \begin{equation}
    \label{eq:ttilde}
      \widetilde T=\{(x^{**},x^*)\in X^{**}\times X^*\;|\;
  (x^*,x^{**})\in  \left(T^{-1}\right)^{\mu,X^*\times X^{**}}\}.
  \end{equation}
  Note that $\rot T=T^{-1}$ and $\rot\,\widetilde T=
  \left(T^{-1}\right)^{\mu,X^*\times X^{**}}$. Therefore
  \begin{equation}
    \label{eq:th.b}
    (\mathcal{S}_T)^*\in\mathcal{F}_{\rot \;\widetilde T}\;.
  \end{equation}
  Moreover, as $\rot$ is a bijection
  which preserve the duality product,
  we conclude that $\widetilde T$ is the
  unique maximal monotone extension of $T$ in $X^{**}\times X^*$.

  Since $T\subset \widetilde T$,
  \begin{align*}
    \varphi_{\rot \widetilde T}(x^*,x^{**})&=\sup_{(y^*,y^{**})\in \rot\widetilde T}
    \inner{x^*}{y^{**}}+\inner{y^*}{x^{**}}-\inner{y^*}{y^{**}}\\
    &=\sup_{(y^{**},y^*)\in \widetilde T}
    \inner{x^*}{y^{**}}+\inner{y^*}{x^{**}}-\inner{y^*}{y^{**}}\\
    &\geq \sup_{(y,y^*)\in T}
    \inner{y}{x^*}+\inner{y^*}{x^{**}}-\inner{y}{y^*}\;\;=(\pi+\delta_T)^*(x^*,x^{**}).
  \end{align*}
  Combining the above equation with the second equality in
  \eqref{eq:ft.tm2} we conclude that $\varphi_{\rot \widetilde T}\geq
  (\mathcal{S}_T)^*$.  Using also the fact that $ \varphi_{\rot\widetilde T}$ is
  minimal in $\mathcal{F}_{\rot \;\widetilde T}$ and \eqref{eq:th.b} we
  obtain  $ \varphi_{\rot\widetilde T}= (\mathcal{S}_T)^*$.

  By Theorem~\ref{th:bs},
  $\varphi_T(x,x^*)=(\mathcal{S}_{\,T})^*(x^*,x)$. Therefore,
  \begin{align*}
    (\varphi_T)^*(x^*,x^{**})
    &= \sup_{(y,y^*)\in X\times X^*}
    \inner{y}{x^*}+\inner{y^*}{x^{**}}
    -\varphi_T(y,y^*) \\
    &= \sup_{(y,y^*)\in X\times X^*}
    \inner{y}{x^*}+\inner{y^*}{x^{**}}
    -(\mathcal{S}_{\,T})^* (y^*,y) \\
    &\leq  \sup_{(y^{**},y^*)\in X^{**}\times X^*}
    \inner{y^{**}}{x^*}+\inner{y^*}{x^{**}}
    -(\mathcal{S}_{\,T})^* ({y^*},y^{**})\\
    &= (\mathcal{S}_{\,T})^{**}(x^{**},x^*).
  \end{align*}
  Take $h\in \mathcal{F}_T$. By Theorem~\ref{th:bs}  $\varphi_T\leq h\leq s_{\,T}$.  Using also the fact that conjugation
  reverts the order, the above equation and assumption~\eqref{eq:cond.as}
 we conclude that, for any
  $(x^*,x^{**})$,
  \begin{align}
    \label{eq:ufa.1}
    \inner{x^*}{x^{**}}\leq \mathcal{S}_{\,T}^*(x^*,x^{**})\leq h^*(x^*,x^{**})\leq \varphi_T^*(x^*,x^{**})\leq
  (\mathcal{S}_{\,T})^{**}(x^{**},x^*).
  \end{align}
  Define $g:X^*\times X^{**}$ as $g=\mathcal{\rot}_{X^*\times X^{**}}
  \,(\mathcal{S}_T)^*$, that is,
  \[ g(x^*,x^{**})=((\mathcal{S}_T)^*)^*(x^{**},x^*).\]
  Using \eqref{eq:th.b} and Theorem~\ref{th:bs.2} we conclude that
  $g\in \mathcal{F}_{\rot\widetilde T}$. Therefore, 
  using again the maximal monotonicity of  $\rot\widetilde T$
  in $X^*\times X^{**}$, we have
  \[  g(x^*,x^{**})=(\mathcal{S}_{\,T})^{**}(x^{**},x^*)=\inner{x^*}{x^{**}},
  \quad \forall (x^*,x^{**})\in \rot\widetilde T.
  \]
  Combining the above equations with \eqref{eq:ufa.1} we conclude that
  $h^*$ majorizes the duality product in $X^*\times X^{**}$ and
  coincides with it in $\rot \widetilde T$. As $h^*$ is also convex and
  closed, we have $h^* \in\mathcal{F}_{\rot\;\widetilde{T}}$.
  
  The fact that $T$ satisfies the restricted Br\o ndsted-Rockafellar property
  follows from the assumption on $\mathcal{S}_{\,T}$ and \cite[Theorem 4.2]{MASvJCA08}.
\end{proof}

\begin{proof}[Proof of Theorem~\ref{converse}]
  Suppose there exist only one $\widetilde T\subset X^{**}\times X^*$ 
   maximal monotone extension of $T$ to 
   $X^{**}\times X^{*}$ and that
   \begin{equation}
     \label{eq:contra}
      \mathcal{S}_T^*(x_0^*,x_0^{**})<\inner{x_0^*}{x_0^{**}}.
   \end{equation}
   As $\rot$ is a bijection that preserves the duality product and
   $\rot T=T^{-1}$, we conclude that $\rot\widetilde T$ is the unique maximal
   monotone extension of $T^{-1}$ to $X^*\times X^{**}$. Using now Theorem~\ref{th:ls-02}, Theorem~\ref{th:legaz-svaiter} and \eqref{eq:ft.tm2} we obtain
   \begin{align}
     \nonumber
     \rot\widetilde T&=\left(T^{-1}\right)^{\mu,X^*\times X^{**}}\\
     \nonumber
     &=\{(x^*,x^{**})\in X^*\times X^{**}\,|\,
     \varphi_{T^{-1}, X^*\times X^{**}}(x^*,x^{**})\leq \inner{x^*}{x^{**}}\}\\
     \label{abc}
    &=\{(x^*,x^{**})\in X^*\times X^{**}\,|\,
     \mathcal{S}_T^*(x^*,x^{**})\leq \inner{x^*}{x^{**}}\}.
   \end{align}

   Suppose that 
   \begin{equation}
     \label{eq:a2}
      (\mathcal{S}_T)^*(x^*,x^{**})<\infty.
   \end{equation}
   Define, for $t\in\R$,
   \[ p(t):= (x_0^{*},x_0^{**})+t(x^*-x_0^{*},x^{**}-x_0^{**})=
    (1-t)(x_0^{*},x_0^{**})+t(x^*,x^{**}).
   \]
   As $(\mathcal{S}_T)^*$ is convex,
   \begin{align*}
      (\mathcal{S}_T)^*(p(t))-\pi_{X^*\times X^{**}}(p(t))\leq& (1-t)(\mathcal{S}_T)^*(x_0^*,x_0^{**})+t
   (\mathcal{S}_{T})^*(x^*,x^{**})\\
    &\;\;-\pi_{X^*\times X^{**}}(p(t)), &\forall t\in[0,1].
   \end{align*}
   As the duality product is continuous, the limit of the right hand
   side of this inequality, for $t\to 0+$ is
   $
   (\mathcal{S}_T)^*(x_0^*,x_0^{**})-\inner{x_0^*}{x_0^{**}}<0.
   $.
   Combining this fact with \eqref{abc} we conclude that
   for $t\geq 0$ and small enough,
   \[
   (x_0^{*},x_0^{**})+t(x^*-x_0^{*},x^{**}-x_0^{**})\in \rot \widetilde{T}. 
   \]
  Altogether, we proved that
   \begin{equation}
     \label{eq:abc}
     \begin{array}{l}
       (\mathcal{S}_T)^*(x^*,x^{**})<\infty\Rightarrow\exists \bar t>0, \;\forall t\in [0,\bar t]\\
       \qquad\qquad\qquad
        (x_0^{*},x_0^{**})+t(x^*-x_0^{*},x^{**}-x_0^{**})\in \rot \widetilde{T}.
     \end{array}
   \end{equation}
   Now, suppose that 
   \[  (\mathcal{S}_T)^*(x_1^*,x_1^{**})<\infty,\qquad  \mathcal{S}_T^*(x_2^*,x_2^{**})<\infty.
   \]
   Then there exists $t>0$ such that
   \[ (x_0^{*},x_0^{**})+t(x^*-x_0^{*},x^{**}-x_0^{**})\in \rot \widetilde{T},\qquad
   (x_1^{*},x_1^{**})+t(x^*-x_1^{*},x^{**}-x_1^{**})\in \rot \widetilde{T},
   \]
   As $\rot \widetilde T$ is (maximal) monotone, the above points are
   monotone related and
   \[ t\inner{x_1^*-x_2^*}{x_1^{**}-x_2^{**}}\geq 0.\]
   Hence, $\inner{x_1^*-x_2^*}{x_1^{**}-x_2^{**}}\geq 0$. Therefore the set
   \[ W:=\{ (x^*,x^{**})\in X^*\times X^{**}\,|\,
   \mathcal{S}_T(x^*,x^{**})<\infty\},
   \]
   is monotone.  By \eqref{abc}, $\rot \widetilde T\subset W$. Hence
   $W=\rot \widetilde T$.
   As $(\mathcal{S}_T)^*$ is convex, $W$ is also convex.
      Therefore, $\rot \widetilde T$ is convex and maximal monotone. Now, using
   Lemma~\ref{bas}  we conclude that $\rot \widetilde T$ is affine. This also
   implies that $\widetilde T$ is affine linear.  As
   \[ T=\widetilde T\cap X\times X^*,\] we conclude that $T$ is affine
   linear, in contradiction with our
   assumptions. Therefore~\eqref{eq:contra} can not hold.
\end{proof}

\begin{proof}[Proof of Theorem \ref{th:main}]
Define 
\begin{equation}
  \label{eq:th.def.s}
  \begin{array}{l}
     \widetilde{T}:=\{(x^{**},x^{*})\in X^{**}\times X^{*}\,| \,
     \inner{x^*-y^*}{x^{**}-y} \geq 0,\; \forall \; (y,y^*)\in T\}\;.
  \end{array}
\end{equation}
Gossez has proved that condition D guarantee that $\widetilde T$ is
the unique maximal monotone extension of $T$ to $X^{**}\times X^*$.

Using again the fact that $\rot$ is a bijection which preserve the
duality product and that $\rot T=T^{-1}$ we conclude that $\rot\,\widetilde
T$ is the unique maximal monotone extension of $T^{-1}$ to $X^*\times
X^{**}$. 

We claim that for any
  $(z^*,z^{**})\in X^*\times X^{**}$
  \begin{align}
    \label{eq:a}
    \inner{x^*}{z^{**}}+\inner{z^*}{x^{**}}- \inner{z^*}{z^{**}}\leq
    (\pi+\delta_T)^{*}(z^*,z^{**}) \quad \forall (x^*,x^{**})\in
    \rot\widetilde T.
  \end{align}
  To prove this claim, take an arbitrary $(x^*,x^{**})\in \rot\widetilde
  T$. Then, $(x^{**},x^{*})\in \widetilde T$ and using
  Definition~\ref{def:g} conclude that there exists a \emph{bounded}
  net $\{ (x_i,x_i^*)\}_{i\in I}\subset T$ such that
  \[ x_i \stackrel{\sigma(X^{**},X^*)}{\longrightarrow} x^{**},\qquad
  x_i^*\stackrel{\norm{\cdot}}{\longrightarrow} x^*.
  \]
  Then, for any $i\in I$
  \[
  \inner{x_i^*}{z^{**}}+\inner{z^*}{x_i}- \inner{z^*}{z^{**}}\leq
  (\pi+\delta_T)^{*}(z^*,z^{**}) .
  \]
  Taking the limit on the above inequality in $i$ we obtain
  the desired result \eqref{eq:a}.
  Taking, in the inequality of \eqref{eq:a}, the $\sup$ on $
  (x^*,x^{**})\in \rot\widetilde T$ we obtain
  \[
  \varphi_{\rot\widetilde T}\leq (\delta_T+\pi)^*.
  \]
  As $\rot\widetilde T$ is maximal monotone in $X^*\times X^{**}$,
  \[ 
  \varphi_{\rot\widetilde T}(x^*,x^{**})\geq \inner{x^*}{x^{**}},\qquad
  \forall (x^*,x^{**})\in X^*\times X^{**}.\]
   Combining the two above inequalities with \eqref{eq:ft.tm2} we conclude that
    \begin{align*}
    (\mathcal{S}_{\,T})^*(x^*,x^{**})
     &\geq  \varphi_{\rot\widetilde T}(x^*,x^{**})\geq
     \inner{x^*}{x^{**}},\qquad \forall (x^*,x^{**})\in X^*\times X^{**}.
  \end{align*}
  To end the proof, use Theorem~\ref{th:1}.
\end{proof}


\begin{thebibliography}{10}

\bibitem{birk}
Garrett Birkhoff.
\newblock {\em Lattice theory}.
\newblock Third edition. American Mathematical Society Colloquium Publications,
  Vol. XXV. American Mathematical Society, Providence, R.I., 1967.

\bibitem{BuIuSv}
R.~S. Burachik, A.~N. Iusem, and B.~F. Svaiter.
\newblock Enlargement of monotone operators with applications to variational
  inequalities.
\newblock {\em Set-Valued Anal.}, 5(2):159--180, 1997.

\bibitem{BuSvteb}
R.~S. Burachik and B.~F. Svaiter.
\newblock {$\epsilon$}-enlargements of maximal monotone operators in {B}anach
  spaces.
\newblock {\em Set-Valued Anal.}, 7(2):117--132, 1999.

\bibitem{BuSvSet02}
R.~S. Burachik and B.~F. Svaiter.
\newblock Maximal monotone operators, convex functions and a special family of
  enlargements.
\newblock {\em Set-Valued Anal.}, 10(4):297--316, 2002.

\bibitem{BuSvProc03}
R.~S. Burachik and B.~F. Svaiter.
\newblock Maximal monotonicity, conjugation and the duality product.
\newblock {\em Proc. Amer. Math. Soc.}, 131(8):2379--2383 (electronic), 2003.

\bibitem{BuSaSv1}
Regina~S. Burachik, Claudia~A. Sagastiz{\'a}bal, and B.~F. Svaiter.
\newblock {$\epsilon$}-enlargements of maximal monotone operators: theory and
  applications.
\newblock In {\em Reformulation: nonsmooth, piecewise smooth, semismooth and
  smoothing methods (Lausanne, 1997)}, volume~22 of {\em Appl. Optim.}, pages
  25--43. Kluwer Acad. Publ., Dordrecht, 1999.

\bibitem{Correa1}
Rafael Correa and Pedro Gajardo.
\newblock Eigenvalues of set-valued operators in {B}anach spaces.
\newblock {\em Set-Valued Anal.}, 13(1):1--19, 2005.

\bibitem{Fitz88}
S.~Fitzpatrick.
\newblock Representing monotone operators by convex functions.
\newblock In {\em Workshop/Miniconference on Functional Analysis and
  Optimization (Canberra, 1988)}, volume~20 of {\em Proc. Centre Math. Anal.
  Austral. Nat. Univ.}, pages 59--65. Austral. Nat. Univ., Canberra, 1988.

\bibitem{GosProc71}
J.-P. Gossez.
\newblock Op\'erateurs monotones non lin\'eaires dans les espaces de {B}anach
  non r\'eflexifs.
\newblock {\em J. Math. Anal. Appl.}, 34:371--395, 1971.

\bibitem{GosProc72}
J.-P. Gossez.
\newblock On the range of a coercive maximal monotone operator in a
  nonreflexive {B}anach space.
\newblock {\em Proc. Amer. Math. Soc.}, 35:88--92, 1972.

\bibitem{GosConvProc76}
J.-P. Gossez.
\newblock On a convexity property of the range of a maximal monotone operator.
\newblock {\em Proc. Amer. Math. Soc.}, 55(2):359--360, 1976.

\bibitem{GosExProc76}
J.-P. Gossez.
\newblock On the extensions to the bidual of a maximal monotone operator.
\newblock {\em Proc. Amer. Math. Soc.}, 62(1):67--71 (1977), 1976.

\bibitem{Str1}
T.~T. Hue and J.~J. Strodiot.
\newblock Convergence analysis of a relaxed extragradient-proximal point
  algorithm application to variational inequalities.
\newblock {\em Optimization}, 54(2):191--213, 2005.

\bibitem{K1}
A.~Kaplan and R.~Tichatschke.
\newblock Bregman-like functions and proximal methods for variational problems
  with nonlinear constraints.
\newblock {\em Optimization}, 56(1-2):253--265, 2007.

\bibitem{Moud01}
Paul-Emile Maing{\'e} and Abdellatif Moudafi.
\newblock A proximal method for maximal monotone operators via discretization
  of a first order dissipative dynamical system.
\newblock {\em J. Convex Anal.}, 14(4):869--878, 2007.

\bibitem{MASvJCA08}
M.~Marques~Alves and B.F. Svaiter.
\newblock Br\o nsted-{R}ockafellar property and maximality of monotone
  operators representable by convex functions in non-reflexive {B}anach spaces.
\newblock {\em Journal of Convex Analysis}, (15), 2008.
\newblock To appear.

\bibitem{legaz-svaiter-01}
J.-E. Mart{\'{\i}}nez-Legaz and B.~F. Svaiter.
\newblock Monotone operators representable by l.s.c.\ convex functions.
\newblock {\em Set-Valued Anal.}, 13(1):21--46, 2005.

\bibitem{SolodovSv1}
M.~V. Solodov and B.~F. Svaiter.
\newblock A hybrid approximate extragradient-proximal point algorithm using the
  enlargement of a maximal monotone operator.
\newblock {\em Set-Valued Anal.}, 7(4):323--345, 1999.

\bibitem{SolodovSv2}
M.~V. Solodov and B.~F. Svaiter.
\newblock A unified framework for some inexact proximal point algorithms.
\newblock {\em Numer. Funct. Anal. Optim.}, 22(7-8):1013--1035, 2001.

\bibitem{SvSet00}
B.~F. Svaiter.
\newblock A family of enlargements of maximal monotone operators.
\newblock {\em Set-Valued Anal.}, 8(4):311--328, 2000.

\bibitem{SvFixProc03}
B.~F. Svaiter.
\newblock Fixed points in the family of convex representations of a maximal
  monotone operator.
\newblock {\em Proc. Amer. Math. Soc.}, 131(12):3851--3859 (electronic), 2003.

\end{thebibliography}
\end{document}